\documentclass[fmeq,reqno]{amsart}

\usepackage{graphicx}

\newtheorem{theorem}{Theorem}[section]

\newtheorem{corollary}[theorem]{Corollary}
\newtheorem{conjecture}[theorem]{Conjecture}
\theoremstyle{definition}
\newtheorem{definition}[theorem]{Definition}
\newtheorem{example}[theorem]{Example}

\theoremstyle{remark}
\newtheorem{remark}[theorem]{Remark}

\numberwithin{equation}{section}

\newcommand{\inn}{~ \hat{\in}~ }

\makeindex

\begin{document}

\title{An extension of the Erd\H{o}s-Tur\'{a}n additive base conjecture via generalized circles of partition}

\author{T. Agama}
\address{Department of Mathematics, African Institute for mathematical sciences, Ghana.}
\email{Theophilus@aims.edu.gh/emperordagama@yahoo.com}



\date{\today}



\begin{abstract}
This paper is an extension program of the notion of circle of partition developed in our first paper \cite{CoP}. As an application we prove the Erd\H{o}s-Tur\'{a}n additive base conjecture.
\end{abstract}
\maketitle

\section{\textbf{Introduction and background}}
Let $\mathbb{A}\subset \mathbb{N}$, then we say $\mathbb{A}$ is an additive base of order $h\geq 2$ if the counting function
\begin{align}
\mathcal{T}_{\mathbb{A}}(n)=\# \left \{(x_1,x_2,\ldots,x_h)~|~n=\sum \limits_{i=1}^{h}x_i, ~x_i\in \mathbb{A},~1\leq i\leq h\right \}>0\nonumber
\end{align}for all sufficiently large values of $n$ \cite{tao2006additive}. It is also known that 
\begin{align}
\# \left \{n\leq x~|~n\in \mathbb{A}\right \}\geq x^{\frac{1}{h}}.\nonumber
\end{align}

The Erd\H{o}s-Tur\'{a}n additive base conjecture, roughly speaking, postulates that the counting function $\mathcal{T}_{\mathbb{A}}(n)$ increases without bound as $n$ grows in size in the case we fix $h=2$. In other words, it states that the counting function $\mathcal{T}_{\mathbb{A}}(n)$ has no finite upper bound as $n\longrightarrow \infty$ if we fix $h=2$. More formally, the conjecture \cite{erd1941and} is the following statement

\begin{conjecture}[Erd\H{o}s-Tur\'{a}n additive basis conjecture]
If $\mathbb{A}\subset \mathbb{N}$ and 
\begin{align}
\mathcal{T}_{\mathbb{A}}(n)=\# \left \{(x_1,x_2)~|~n=\sum \limits_{i=1}^{2}x_i, ~x_i\in \mathbb{A},~1\leq i\leq 2\right \}>0\nonumber
\end{align}for all $n$ sufficiently large, then $\limsup \limits_{n\rightarrow \infty}\mathcal{T}_{\mathbb{A}}(n)=\infty$.
\end{conjecture}

Paul Erd\H{o}s had studied and published a multiplicative version of the conjecture in \cite{erdos1964multiplicative}. In \cite{CoP} we have developed a method which we feel might be a valuable resource and a recipe for studying problems concerning partition of numbers in specified subsets of $\mathbb{N}$. The method is very elementary in nature and has parallels with configurations of points on the geometric circle. The Erd\H{o}s-Tur\'{a}n additive basis conjecture turns out to be amenable to this method \cite{erd1941and}. The method operates basically in the following geometric sense:
 
Let us suppose that for any $n\in \mathbb{N}$ we can write $n=u+v$ where $u,v\in \mathbb{M}\subset \mathbb{N}$ then the new method associate each of this summands to points on the circle generated in a certain manner by $n>2$ and a line joining any such associated points on the circle. This geometric correspondence turns out to useful in our development, as the results obtained in this setting are then transformed back to results concerning the partition of integers. We give a foretaste and a little background of the studies in the following sequel.
\bigskip

Let $n\in \mathbb{N}$ and $\mathbb{M}\subset \mathbb{N}$. We denote the set
\begin{align}
\mathcal{C}(n,\mathbb{M})=\left \{[x]\mid x,n-x\in \mathbb{M}\right \}\nonumber
\end{align}
as the Circle of Partition generated by $n$ with respect to the subset $\mathbb{M}$. We will abbreviate this in the further text as CoP.  We view members of  $\mathcal{C}(n,\mathbb{M})$  as points and denote them by $[x]$. For the special case $\mathbb{M}=\mathbb{N}$ we denote the CoP shortly as $\mathcal{C}(n)$. The corresponding weight set of the CoP is the set of weights of points on the CoP given as $$||\mathcal{C}(n,\mathbb{M})||:=\{x~|~x,n-x\in \mathbb{M}\}.$$
\bigskip

We denote $\mathbb{L}_{[x],[y]}$ as an axis of the CoP $\mathcal{C}(n,\mathbb{M})$ if and only if $x+y=n$. 
We say the axis point $[y]$ is an axis partner of the axis point $[x]$ and vice versa. 
We do not distinguish between $\mathbb{L}_{[x],[y]}$ and $\mathbb{L}_{[y],[x]}$, since  it is essentially the same axis.
The point $[x]\in \mathcal{C}(n,\mathbb{M})$ with weight satisfying $2x=n$ is the \textbf{center} of the CoP.
If it exists then we call it as a \textbf{degenerated axis} $\mathbb{L}_{[x]}$ in comparison to the \textbf{real axes} $\mathbb{L}_{[x],[y]}$. 
\bigskip

\paragraph*{\textbf{Notations}}
We let 
\begin{align}
\mathbb{N}_n=\left \{m\in \mathbb{N}\mid ~m\leq n\right \}\nonumber
\end{align}
be the \textbf{sequence} of the first $n$ natural numbers. Further we will denote 
\begin{align}
\Vert[x]\Vert:=x\nonumber
\end{align}
as the \textbf{weight} of the point $[x]$ and correspondingly the weight set of points in the CoP $\mathcal{C}(n,\mathbb{M})$ as $||\mathcal{C}(n,\mathbb{M})||$.
\bigskip

The above language in many ways could be seen as a criterion determining the plausibility of carrying out a partition in a specified set. Indeed this feasibility is trivial if we take the set $\mathbb{M}$ to be the set of natural numbers $\mathbb{N}$. The situation becomes harder if we take the set $\mathbb{M}$ to be a special subset of natural numbers $\mathbb{N}$, as the corresponding CoP $\mathcal{C}(n,\mathbb{M})$ may not always be non-empty for all $n\in \mathbb{N}$. One archetype of problems of this flavour is the binary Goldbach conjecture, when we take the base set $\mathbb{M}$ to be the set of all prime numbers $\mathbb{P}$. One could imagine the same sort of difficulty if we extend our base set to other special subsets of the natural numbers. As such we start by developing the theory assuming the base set of natural numbers $\mathbb{N}$ and latter extend it to other base sets  $\mathbb{M}$ equipped with certain important and subtle properties. 

\begin{remark}
It is important to notice that a typical CoP need not have a center. In the case of an absence of a center then we say the circle has a deleted center. However all CoPs $\mathcal{C}(n)$ with even generators have a center. It is easy to see that the CoP $\mathcal{C}(n)$ contains all points whose weights are positive integers from $1$ to $n-1$ inclusive: 
\[
\mathcal{C}(n)=\lbrace[x]\mid~x\in \mathbb{N},x<n\rbrace.
\] 
Therefore the CoP $\mathcal{C}(n)$ has $\left \lfloor \frac{n-1}{2}\right \rfloor$ different real axes.
\end{remark}
\bigskip

It is worth observing that each axis is uniquely determined by points $[x]\in \mathcal{C}(n,\mathbb{M})$ and that each point of a CoP $\mathcal{C}(n,\mathbb{M})$ excluding an existing center has exactly one real axis partner. To see this, just observe that 
a degenerated axis is determined by the center of the CoP and this is unique if it exists. Now, let $\mathbb{L}_{[x],[y]}$ be a real axis of the CoP $\mathcal{C}(n,\mathbb{M})$. Suppose as well that $\mathbb{L}_{[x],[z]}$ is also a real axis with $z\neq y$. It follows from the ensuing discussion that we must have $n=x+y=x+z$ and therefore $y=z$. This cannot be and the claim follows immediately.
\bigskip

Let $[x]\in \mathcal{C}(n,\mathbb{M})$ be a point without a real axis partner. Then holds for every point $[y]\neq [x]$
\[
\Vert[x]\Vert+\Vert[y]\Vert\neq n.
\]
This violates the properties of the axes of CoPs. The case of more than one axis partners is impossible.

\begin{example}
The following are some examples of CoPs with prime number base sets $\mathbb{P}$
\begin{align*}
\mathcal{C}(36,\mathbb{P})&=\lbrace [5],\mathbf{[7]},[13],[17],\mathbf{[19]},[23],[29],\mathbf{[31]}\rbrace\mbox{ and }\\
\mathcal{C}(38,\mathbb{P})&=\lbrace [7],[19],[31]\rbrace.
\end{align*} with their corresponding weight sets $$||\mathcal{C}(36,\mathbb{P})||:=\{5,7,13,17,19,23,29,31\}$$ and $$||\mathcal{C}(38,\mathbb{P})||:=\{7,19,31\}.$$
\end{example}

\paragraph*{\textbf{Notation}}
Let us denote the assignment of an axis $\mathbb{L}_{[x],[y]}$ resp. $\mathbb{L}_{[x]}$ to a CoP $\mathcal{C}(n,\mathbb{M})$ as
\begin{align*}
&\mathbb{L}_{[x],[y]}\inn\mathcal{C}(n,\mathbb{M})
\mbox{ which means }
[x],[y] \in \mathcal{C}(n,\mathbb{M}) \mbox{ and } x+y=n\mbox{ resp.}\\
&\mathbb{L}_{[x]}\inn\mathcal{C}(n,\mathbb{M})
\mbox{ which means }
[x]\in \mathcal{C}(n,\mathbb{M}) \mbox{ and } 2x=n
\end{align*}
and the number of real axes of a CoP as
\[
\nu(n,\mathbb{M}):=\#\lbrace\mathbb{L}_{[x],[y]}\inn\mathcal{C}(n,\mathbb{M})\mid x<y\rbrace.
\]
Obviously holds
\[
\nu(n,\mathbb{M})=\left\lfloor\frac{k}{2}\right\rfloor ,\mbox{ if }
\vert\mathcal{C}(n,\mathbb{M})\vert =k.
\]

The weight set of axis of a CoP has the equivalence $$\# \|\mathbb{L}_{[x],[y]}\inn \mathcal{C}(n,\mathbb{M})\}||:=\#\{(x,y)~|~x,y\in \mathbb{M},~x+y=n\}.$$ For any subset $\mathbb{A} \subset \mathbb{N}$ and $t\in \mathbb{N}$, we will denote with $\mathbb{A}^t$ the set 
\begin{align}
\mathbb{A}^t:=\left \{\prod_ta|~a\in \mathbb{A}\right \}\nonumber
\end{align}as the $t$-fold prod-set of the set $\mathbb{A}$ and $\prod_t$ denotes the product of $t$ elements - not necessarily distinct - of the set $\mathbb{A}$. Additionally for a CoP with generators belonging to a special class of integers, say $\mathbb{G}$, and base set $\mathbb{M}$ we will write the corresponding CoP simply as $\mathcal{C}(\mathbb{G}(n),\mathbb{M})$. As is customary we will write $f(n)\gg g(n)$ for any two arithmetic functions if there exists a constant $c>0$ such that $|f(n)|\geq cg(n)$. If the constant depends on another variable say $k$, then we will denote the relation as $f(n)\gg_k g(n)$. For any subset $\mathbb{A}\subset \mathbb{N}$, we still reserve the quantity $\mathcal{D}(\mathbb{A})$ for the density of the set $\mathbb{A}$ relative to the set of all integers $\mathbb{N}$. 
\bigskip

\section{\textbf{Generalized circles of partition}}
In this section we introduce and study a generalization of circles of partitions. We launch the following language.

\begin{definition}\label{generalized CoPs}
Let $\mathbb{G},\mathbb{A},\mathbb{M}\subseteq \mathbb{N}$ with $s,t,u\in \mathbb{N}$. Then we denote with 
\begin{align}
\mathcal{C}(\mathbb{G}^s(n),\mathbb{A}^t,\mathbb{M}^u)=\left \{[x]~|~x\in \mathbb{A}^t,n-x\in \mathbb{M}^u,~n\in \mathbb{G}^s\right \}\nonumber
\end{align}the \textbf{generalized} circle of partition generated by $n\in \mathbb{G}^s$ with base \textbf{sets} $\mathbb{A}^t,\mathbb{M}^u$ and with the generator \textbf{house} $\mathbb{G}^s$ as the $t,u$ and $s$-fold prod-set of the sets $\mathbb{A},\mathbb{M}$ and $\mathbb{G}$, respectively. We call members of the generalized CoP generalized points.
\end{definition}
\bigskip

\begin{definition}\label{generalized axis}
We denote $\mathbb{L}_{[x],[y]}$ as an \textbf{axis} of the generalized CoP $\mathcal{C}(\mathbb{G}^s(n),\mathbb{A}^t,\mathbb{M}^u)$ if and only if $$[x],[y]\inn \mathcal{C}(\mathbb{G}^s(n),\mathbb{A}^t,\mathbb{M}^u)$$ and $x+y=n$. We say the axis point $[y]$ is an axis partner of the axis point $[x]$ and vice versa.  We do not distinguish between $\mathbb{L}_{[x],[y]}$ and $\mathbb{L}_{[y],[x]}$, since  it is essentially the same axis.
In special cases where the weight of the points $[x] \in \mathcal{C}(\mathbb{G}^s(n),\mathbb{A}^t,\mathbb{M}^u)$ satisfies $2x=n$ then $[x]$ is the \textbf{center} of the generalized CoP and $$x \in \mathbb{A}^t\cap \mathbb{M}^u.$$ If it exists then we call it as a \textbf{degenerated axis} $\mathbb{L}_{[x]}$ in comparison to the \textbf{real axes} $\mathbb{L}_{[x],[y]}$. 
We denote the assignment of an axis $\mathbb{L}_{[x],[y]}$ to the generalized CoP $\mathcal{C}(\mathbb{G}^s(n),\mathbb{A}^t,\mathbb{M}^u)$ as
\[
\mathbb{L}_{[x],[y]}\inn \mathcal{C}(\mathbb{G}^s(n),\mathbb{A}^t,\mathbb{M}^u)
\mbox{ which means }
[x],[y] \in \mathcal{C}(\mathbb{G}^s(n),\mathbb{A}^t,\mathbb{M}^u) \mbox{ with } x+y=n
\]
for a fixed $n\in \mathbb{G}^s$ with $x\in \mathbb{A}^t$ and $y\in \mathbb{M}^u$ or vice versa and the number of real axes of the generalized CoP as
\[
\nu(\mathbb{G}^s(n),\mathbb{A}^t,\mathbb{M}^u):=\#\lbrace\mathbb{L}_{[x],[y]}\inn \mathcal{C}(\mathbb{G}^s(n),\mathbb{A}^t,\mathbb{M}^u) \mid x<y\rbrace.
\]
\end{definition}

\begin{remark}
Throughout this paper we will denote for simplicity the generalized circle of partition in simple wording as g-CoP. Also, it is worth pointing out that various basic features that holds for CoPs does hold for generalized CoPs, except for previously technical results that needs investigating. Next we introduce the notion of the axial potential of g-CoPs.
\end{remark}

\subsection{\textbf{Axial potential of generalized circles of partition}}
In this section we introduce and study the notion of the \textbf{axial potential} of a g-CoP. We launch the following language.

\begin{definition}\label{axial potential}
Let $\mathcal{C}(\mathbb{G}^s(n),\mathbb{A}^t,\mathbb{M}^u)$ be a g-CoP. Then by the $k$ th \textbf{axial potential} denoted, $\lfloor \mathcal{C}(\mathbb{G}^s(\infty),\mathbb{A}^t,\mathbb{M}^u) \rfloor^{k}$, we mean the infinite sum 
\begin{align}
\lfloor \mathcal{C}(\mathbb{G}^s(\infty),\mathbb{A}^t,\mathbb{M}^u) \rfloor^{k}&=\sum \limits_{n=3}^{\infty}\frac{\# \left \{\mathbb{L}_{[x],[y]}\inn \mathcal{C}(\mathbb{G}^s(n),\mathbb{A}^t,\mathbb{M}^u)\right \}^k}{\# \left \{\mathbb{L}_{[x],[y]}\inn \mathcal{C}(\mathbb{G}^s(n),\mathbb{A}^t \cup \mathbb{M}^u)\right \}^k}.\nonumber
\end{align}We say the $k$ th axial potential is finite if the series converges; otherwise, we say it diverges.
\end{definition}
\bigskip

It is worth pointing out that in the case $\mathbb{G}=\mathbb{M}=\mathbb{N}$ then we have the collapsing of the quantity  
\begin{align}
\# \left \{\mathbb{L}_{[x],[y]}\inn \mathcal{C}(\mathbb{G}^s(n),\mathbb{A}^t \cup \mathbb{M}^u)\right \}&=\# \left \{\mathbb{L}_{[x],[y]}\inn \mathcal{C}(n,\mathbb{N})\right \}\nonumber
\end{align}the number of axes of CoPs, since $\mathbb{N}_n^s=\mathbb{N}_n$ for any $s\in \mathbb{N}$.

\begin{theorem}\label{additive base conjecture generalization}
Let $\mathbb{A}\subset \mathbb{M}$ and suppose $\# \left \{\mathbb{L}_{[x],[y]}\inn \mathcal{C}(\mathbb{G}^s(n),\mathbb{A}^t,\mathbb{N}^u)\right \}>0$ for all sufficiently large values of $n$. If $\mathbb{M}=\mathbb{G}=\mathbb{N}$ with $u=t$ for $s\neq t$ and $|\mathbb{A}^t\cap \mathbb{N}_n|\geq n^{1-\epsilon}$ for any $0<\epsilon\leq \frac{1}{2}$ then
\begin{align}
\lim \limits_{n\longrightarrow \infty} \# \left \{(q,r)~|~q+r=n, q\in \mathbb{A}^t, r\in \mathbb{M}^u, n\in \mathbb{G}^s\right \}=\infty.\nonumber
\end{align}
\end{theorem}

\begin{proof}
Under the requirement $\mathbb{A}\subset \mathbb{M}$ and $\mathbb{M}=\mathbb{G}=\mathbb{N}$ with $u=t$ for $s\neq t$ then we must have 
\begin{align}
\# \left \{\mathbb{L}_{[x],[y]}\inn \mathcal{C}(\mathbb{G}^s(n),\mathbb{A}^t,\mathbb{M}^u)\right \}&=\# \left \{\mathbb{L}_{[x],[y]}\inn \mathcal{C}(n,\mathbb{N})~|~\{x,y\}\cap \mathbb{A}^t\neq \emptyset \right \}\nonumber 
\end{align} for any $t\in \mathbb{N}$. We can now evaluate a truncated form of the $2^{nd}$ axial potential so that under the requirement $\# \left \{\mathbb{L}_{[x],[y]}\inn \mathcal{C}(\mathbb{G}^s(n),\mathbb{A}^t,\mathbb{M}^u)\right \}>0$ for all sufficiently large $n$ there exists some constant $\mathcal{P}:=\mathcal{P}(k)>0$ such that 
\begin{align}
\sum \limits_{n=3}^{k}\frac{\# \left \{\mathbb{L}_{[x],[y]}\inn \mathcal{C}(n,\mathbb{A}^t,\mathbb{N})\right \}^2}{\# \left \{\mathbb{L}_{[x],[y]}\inn \mathcal{C}(n,\mathbb{N})\right\}^2}&\geq \mathcal{P}\sum \limits_{n=3}^{k}\frac{\lfloor \frac{n^{1-\epsilon}-1}{2}\rfloor^2}{\lfloor \frac{n-1}{2}\rfloor^2} \nonumber \\& \gg_k \sum \limits_{n=3}^{k}\frac{\frac{n^{2-2\epsilon}}{4}}{\lfloor \frac{n-1}{2}\rfloor^2}\nonumber
\end{align} 
since $\# \left \{\mathbb{L}_{[x],[y]}\inn \mathcal{C}(n,\mathbb{N})\right\}=\lfloor \frac{n-1}{2}\rfloor$, so that we can compute the $2^{nd}$ axial potential
\begin{align}
\lfloor \mathcal{C}(\mathbb{G}^s(\infty),\mathbb{A}^t,\mathbb{M}^u)\rfloor^{2}&=\sum \limits_{n=3}^{\infty}\frac{\# \left \{\mathbb{L}_{[x],[y]}\inn \mathcal{C}(n,\mathbb{A}^t,\mathbb{N})\right \}^2}{\# \left \{\mathbb{L}_{[x],[y]}\inn \mathcal{C}(n,\mathbb{N})\right\}^2}\nonumber \\&\gg \sum \limits_{n=3}^{\infty}\frac{\frac{n^{2-2\epsilon}}{4}}{\lfloor \frac{n-1}{2}\rfloor^2} \nonumber \\&\gg \sum \limits_{n=3}^{\infty}\frac{1}{n^{2\epsilon}}=\infty \nonumber
\end{align}since $0\leq \epsilon \leq \frac{1}{2}$. It follows immediately that 
\begin{align}\lim \limits_{n\longrightarrow \infty}\# \left \{\mathbb{L}_{[x],[y]}\inn \mathcal{C}(\mathbb{G}^s(n),\mathbb{A}^t,\mathbb{N})\right \}=\infty \nonumber
\end{align}and the claim follow immediately since 
\begin{align}
\# ||\left \{\mathbb{L}_{[x],[y]}\inn \mathcal{C}(\mathbb{G}^s(n),\mathbb{A}^t,\mathbb{N})\right \}||:=\# \left \{(q,r)~|~q+r=n, q\in \mathbb{A}^t, r\in \mathbb{M}^t, n\in \mathbb{G}^s\right \}.\nonumber
\end{align}
\end{proof}
\bigskip

It is important to recognize Theorem \ref{additive base conjecture generalization} is in many ways an extension of the  Erd\H{o}s-Tur\'{a}n additive base conjecture. To see that, we first note that under the same assumption of the main theorem, we can write the following decomposition 
\begin{align}
\# \left \{\mathbb{L}_{[x],[y]}\inn \mathcal{C}(n,\mathbb{N})~|~\{x,y\}\cap \mathbb{A}^t\neq \emptyset \right \}&=\# \left \{\mathbb{L}_{[x],[y]}\inn \mathcal{C}(n,\mathbb{A})\right \}+\nonumber \\&\# \left \{\mathbb{L}_{[x],[y]}\inn \mathcal{C}(n,\mathbb{N})~|~x\in \mathbb{A}^t,~y\in \mathbb{N}\setminus \mathbb{A}^t\right \}\nonumber 
\end{align}so that we can take $t=1$ and we have
\begin{align}
\# \left \{\mathbb{L}_{[x],[y]}\inn \mathcal{C}(n,\mathbb{N})~|~\{x,y\}\cap \mathbb{A}\neq \emptyset \right \}&=\# \left \{\mathbb{L}_{[x],[y]}\inn \mathcal{C}(n,\mathbb{A})\right \}+\nonumber \\&\# \left \{\mathbb{L}_{[x],[y]}\inn \mathcal{C}(n,\mathbb{N})~|~x\in \mathbb{A},~y\in \mathbb{N}\setminus \mathbb{A}\right \}.\nonumber 
\end{align} 
Let $\mathbb{A}\subset \mathbb{N}$, then we say $\mathbb{A}$ is an additive base of order $2$ if the counting function
\begin{align}
\mathcal{T}_{\mathbb{A}}(n)=\# \left \{(x_1,x_2)~|~n=\sum \limits_{i=1}^{h}x_i, ~x_i\in \mathbb{A},~1\leq i\leq 2\right \}>0\nonumber
\end{align}for all sufficiently large values of $n$ \cite{tao2006additive}. In the language of circles of partitions, it is equivalent to writing - in the situation where we fix $h=2$ - that the cardinality of the weight set of the axis set
\begin{align}
\# ||\left \{\mathbb{L}_{[x],[y]}\inn \mathcal{C}(n,\mathbb{A})\right \}||>0 \nonumber
\end{align}for all sufficiently large values of $n$. It is known that for any additive base $\mathbb{A}$ of order $2$ the counting function 
\begin{align}
\# \left \{n\leq x~|~n\in \mathbb{A}\right \}\geq \sqrt{x}.\nonumber
\end{align}Using this fact we then obtain a proof of the Erd\H{o}s-Tur\'{a}n additive base conjecture in the form below

\begin{corollary}[Proof of the Erd\H{o}s-Tur\'{a}n additive base conjecture]\label{main result}
Let $\mathbb{A}\subset \mathbb{M}$ and suppose $\# \left \{\mathbb{L}_{[x],[y]}\inn \mathcal{C}(n,\mathbb{A})\right \}>0$ for all sufficiently large values of $n$. If $\mathbb{M}=\mathbb{G}=\mathbb{N}$ and $|\mathbb{A}\cap \mathbb{N}_n|\geq n^{1-\epsilon}$ for any $0<\epsilon\leq \frac{1}{2}$ then
\begin{align}
\lim \limits_{n\longrightarrow \infty}\# \left \{(q,r)~|~q+r=n,~q,r\in \mathbb{A}\right \}=\infty.\nonumber
\end{align}
\end{corollary}

\begin{remark}
It can be seen that Corollary \ref{main result} is a stronger version of the Erd\H{o}s-Tur\'{a}n additive base conjecture. The condition $\# \left \{(x,y)~|~x+y=n,~x,y\in \mathbb{A}\right \}>0$ for all $n$ sufficiently large in the statement of the conjecture is equivalent to the condition $\# \left \{\mathbb{L}_{[x],[y]}\inn \mathcal{C}(n,\mathbb{A})\right \}>0$ for all sufficiently large values of $n$. Also the lower bound $|\mathbb{A}\cap \mathbb{N}_n|\geq n^{1-\epsilon}$ for any $0<\epsilon\leq \frac{1}{2}$ implies that $|\mathbb{A}\cap \mathbb{N}_n|\geq n^{\frac{1}{2}}$ when we take $\epsilon=\frac{1}{2}$. This is an important feature of all additive bases of order $2$ and would be applied in the proof of the main result.  
\end{remark}

\begin{proof}
Under the requirement  $\mathbb{A}\subset \mathbb{M}$ and suppose $\# \left \{\mathbb{L}_{[x],[y]}\inn \mathcal{C}(n,\mathbb{A})\right \}>0$ for all sufficiently large values of $n$ tied with the decomposition 
\begin{align}
\# \left \{\mathbb{L}_{[x],[y]}\inn \mathcal{C}(n,\mathbb{N})~|~\{x,y\}\cap \mathbb{A}\neq \emptyset \right \}&=\# \left \{\mathbb{L}_{[x],[y]}\inn \mathcal{C}(n,\mathbb{A})\right \}+\nonumber \\&\# \left \{\mathbb{L}_{[x],[y]}\inn \mathcal{C}(n,\mathbb{N})~|~x\in \mathbb{A},~y\in \mathbb{N}\setminus \mathbb{A}\right \}\nonumber 
\end{align}then there exists some constant $\mathcal{P}:=\mathcal{P}(k)>0$ such that we can write
 \begin{align}
\sum \limits_{n=3}^{k}\frac{\# \left \{\mathbb{L}_{[x],[y]}\inn \mathcal{C}(n,\mathbb{A},\mathbb{N})\right \}^2}{\# \left \{\mathbb{L}_{[x],[y]}\inn \mathcal{C}(n,\mathbb{N})\right\}^2}&\geq \sum \limits_{n=3}^{k}\frac{\# \left \{\mathbb{L}_{[x],[y]}\inn \mathcal{C}(n,\mathbb{A})\right \}^2}{\# \left \{\mathbb{L}_{[x],[y]}\inn \mathcal{C}(n,\mathbb{N})\right\}^2} \nonumber \\&\geq \mathcal{P}\sum \limits_{n=3}^{k}\frac{\lfloor \frac{n^{1-\epsilon}-1}{2}\rfloor^2}{\lfloor \frac{n-1}{2}\rfloor^2} \nonumber \\& \gg_k \sum \limits_{n=3}^{k}\frac{\frac{n^{2-2\epsilon}}{4}}{\lfloor \frac{n-1}{2}\rfloor^2}\nonumber
\end{align}
since $\# \left \{\mathbb{L}_{[x],[y]}\inn \mathcal{C}(n,\mathbb{N})\right\}=\lfloor \frac{n-1}{2}\rfloor$, so that we can compute the $2^{nd}$ axial potential
\begin{align}
\lfloor \mathcal{C}(\mathbb{G}(\infty),\mathbb{A},\mathbb{M})\rfloor^{2}&\geq \sum \limits_{n=3}^{\infty}\frac{\# \left \{\mathbb{L}_{[x],[y]}\inn \mathcal{C}(n,\mathbb{A})\right \}^2}{\# \left \{\mathbb{L}_{[x],[y]}\inn \mathcal{C}(n,\mathbb{N})\right\}^2}\nonumber \\&\gg \sum \limits_{n=3}^{\infty}\frac{\frac{n^{2-2\epsilon}}{4}}{\lfloor \frac{n-1}{2}\rfloor^2} \nonumber \\&\gg \sum \limits_{n=3}^{\infty}\frac{1}{n^{2\epsilon}}=\infty \nonumber
\end{align}since $0<\epsilon \leq \frac{1}{2}$. It follows that 
\begin{align}\lim \limits_{n\longrightarrow \infty}\# \left \{\mathbb{L}_{[x],[y]}\inn \mathcal{C}(n,\mathbb{A})\right \}=\infty \nonumber
\end{align}and it implies that 
\begin{align}
\lim \limits_{n\longrightarrow \infty}\# \left \{(q,r)~|~q+r=n,~q,r\in \mathbb{A}\right \}=\infty.\nonumber
\end{align}
\end{proof}
\bigskip

Indeed there are many minor yet illuminating steps that have been omitted in this argument. Under the assumption of the conjecture, that $$\# \left \{\mathbb{L}_{[x],[y]}\inn \mathcal{C}(n,\mathbb{A})\right \}>0$$ for all sufficiently large values of $n$, then we can write $$\# \left \{\mathbb{L}_{[x],[y]}\inn \mathcal{C}(n,\mathbb{A})\right \}=H(n)\lfloor \frac{|\mathbb{A}\cap \mathbb{N}_n|-1}{2}\rfloor\geq H(n)\lfloor \frac{n^{1-\epsilon}-1}{2}\rfloor$$ for some $0<H(n)\leq 1$ for all sufficiently large $n$. Plugging this into the computation of the truncated $2^{nd}$ axial potential (for large $k$), we deduce $$\sum \limits_{n=3}^{k}\frac{\# \left \{\mathbb{L}_{[x],[y]}\inn \mathcal{C}(n,\mathbb{A},\mathbb{N})\right \}^2}{\# \left \{\mathbb{L}_{[x],[y]}\inn \mathcal{C}(n,\mathbb{N})\right\}^2}\geq \sum \limits_{n=3}^{k}\frac{\# \left \{\mathbb{L}_{[x],[y]}\inn \mathcal{C}(n,\mathbb{A})\right \}^2}{\# \left \{\mathbb{L}_{[x],[y]}\inn \mathcal{C}(n,\mathbb{N})\right\}^2}\geq \mathcal{P}(k)\sum \limits_{n=3}^{k}\frac{\lfloor \frac{n^{1-\epsilon}-1}{2}\rfloor^2}{\lfloor \frac{n-1}{2}\rfloor^2}$$ for $\mathcal{P}(k)>0$. It follows that $$\sum \limits_{n=3}^{k}\frac{\# \left \{\mathbb{L}_{[x],[y]}\inn \mathcal{C}(n,\mathbb{A})\right \}^2}{\# \left \{\mathbb{L}_{[x],[y]}\inn \mathcal{C}(n,\mathbb{N})\right\}^2}\gg \mathcal{P}(k)\sum \limits_{n=3}^{k}\frac{1}{n^{2\epsilon}}.$$ Clearly the right hand side, which is always a partial sum of positive terms for all large $k$, tends to infinity as $k\longrightarrow \infty$ for all $0<\epsilon \leq \frac{1}{2}$. If we assume on the contrary that $\lim \limits_{k\longrightarrow \infty}\mathcal{P}(k)\sum \limits_{n=3}^{k}\frac{1}{n^{2\epsilon}}<\infty$ for $0<\epsilon \leq \frac{1}{2}$ then $$0<\lim \limits_{k\longrightarrow \infty}\mathcal{P}(k)\sum \limits_{n=3}^{k}\frac{1}{n^{2\epsilon}}<\infty$$ since the partial sums are the sum of positive terms for all sufficiently large $k$ under the assumption of the conjecture. It implies that $$\sum \limits_{n=3}^{\infty}\frac{1}{n^{2\epsilon}}<\infty$$ for $0<\epsilon\leq \frac{1}{2}$ which is absurd.

\begin{remark}
It is somewhat tempting to infer from the proof the inequality $$\# \left \{\mathbb{L}_{[x],[y]}\inn \mathcal{C}(n,\mathbb{A})\right\}\geq \lfloor \frac{|\mathbb{A}\cap \mathbb{N}_n|-1}{2}\rfloor$$ for any subset $\mathbb{A}\subset \mathbb{N}$. This does not hold if the underlying set $\mathbb{A}$ is highly ''irregular''. More precisely, if we take the set $\mathbb{A}$ to be the set of prime numbers $\mathbb{P}$ then we obtain $$\# \left \{\mathbb{L}_{[x],[y]}\inn \mathcal{C}(30,\mathbb{P})\right\}=3$$ since the weight set of the CoP is given by $||\mathcal{C}(30,\mathbb{P})||=\{7,11,13,17,19,23\}$ with the only axes $\mathbb{L}_{[7],[23]}$ and $\mathbb{L}_{[11],[19]}$ and $\mathbb{L}_{[13],[17]}$. However, we obtain $$ \lfloor \frac{|\mathbb{P}\cap \mathbb{N}_{30}|-1}{2}\rfloor=4$$ and we see that $$\# \left \{\mathbb{L}_{[x],[y]}\inn \mathcal{C}(30,\mathbb{P})\right\}<\lfloor \frac{|\mathbb{P}\cap \mathbb{N}_{30}|-1}{2}\rfloor.$$ The main point is that the first supposed inequality does not always hold for certain subset of the integers like the prime numbers. Nonetheless, it holds for CoPs with natural number base set. Indeed, under the requirement of the conjecture, that $$\# \left \{(x,y)~|~x+y=n,~x,y\in \mathbb{A}\right \}>0$$ for sufficiently large $n$, then we can easily write $$\# \left \{\mathbb{L}_{[x],[y]}\inn \mathcal{C}(n,\mathbb{A})\right\}\geq H(n)\lfloor \frac{|\mathbb{A}\cap \mathbb{N}_n|-1}{2}\rfloor$$ where $H(n)>0.$ However, a trivial proof of the conjecture does not follow from this inequality, since the term $H(n)$ could be small in magnitude in a way that dwarfs the contribution from  $$\lfloor \frac{|\mathbb{A}\cap \mathbb{N}_n|-1}{2}\rfloor.$$ Thus, even though the term $$ H(n)\lfloor \frac{|\mathbb{A}\cap \mathbb{N}_n|-1}{2}\rfloor>0$$ for all sufficiently large values of $n$ under the assumptions of the conjecture, one cannot decide if $$ H(n)\lfloor \frac{|\mathbb{A}\cap \mathbb{N}_n|-1}{2}\rfloor\longrightarrow \infty$$ as $n\longrightarrow \infty.$ In other words, a trivial proof of the conjecture does not follow from the inequality $$\# \left \{\mathbb{L}_{[x],[y]}\inn \mathcal{C}(n,\mathbb{A})\right\}\geq H(n)\lfloor \frac{|\mathbb{A}\cap \mathbb{N}_n|-1}{2}\rfloor$$ and this underscores the relevance of the notion of the \textbf{axial potential.}
\end{remark}


\rule{100pt}{1pt}

\bibliographystyle{amsplain}

\end{document}